# EVERY ORIENTABLE SEIFERT 3-MANIFOLD IS A REAL COMPONENT OF A UNIRULED ALGEBRAIC VARIETY

JOHANNES HUISMAN AND FRÉDÉRIC MANGOLTE

ABSTRACT. We show that any orientable Seifert 3-manifold is diffeomorphic to a connected component of the set of real points of a uniruled real algebraic variety, and prove a conjecture of János Kollár.

### Toute 3-variété de Seifert est une composante réelle d'une variété algébrique uniréglée

RÉSUMÉ. Nous montrons que toute 3-variété de Seifert orientable est difféomorphe à une composante connexe de la partie réelle d'une variété algébrique réelle uniréglée et prouvons une conjecture de János Kollár.

## 1. INTRODUCTION

A smooth compact connected 3-manifold $M$ is a *Seifert manifold* if it admits a smooth fibration $f\colon M \to S$ over a smooth surface $S$, whose fibers are circles, such that $f$ is locally trivial, with respect to the ramified Grothendieck topology on $S$. More precisely, for every point $P$ of $S$, there is a ramified covering $\tilde{U} \to U$ of an open neighborhood $U$ of $P$ such that the fiber product $M \times_S \tilde{U} \to \tilde{U}$ is a locally trivial smooth circle bundle over $\tilde{U}$ (see Section 2 for another—but equivalent—definition).

A smooth, projective and geometrically irreducible real algebraic variety $X$ is called *ruled* if there is a real algebraic variety $Y$ such that $Y \times \mathbb{P}^1$ and $X$ are birational. The variety $X$ is *uniruled* if there is a real algebraic variety $Y$, with $\dim(Y) = \dim(X) - 1$, and a dominant rational map $Y \times \mathbb{P}^1 \dashrightarrow X$. Of course, a ruled real algebraic variety is uniruled, but not conversely.

Let $X$ be a uniruled real algebraic variety of dimension 3 such that $X(\mathbb{R})$ is orientable. János Kollár has proved that each connected component of $X(\mathbb{R})$ belongs to a given list of manifolds, containing the Seifert manifolds [Ko01, Th. 6.6]. He conjectured, conversely, that each orientable Seifert manifold is diffeomorphic to a connected component of the set of real points of a uniruled real algebraic variety [Ko01, Conj. 6.7.2]. In this paper we prove that conjecture.

**Theorem 1.1.** *Every orientable Seifert manifold is diffeomorphic to a real component of a uniruled real algebraic variety.*

The strategy of our proof is the following. Let $M$ be an orientable Seifert manifold, and let $f\colon M \to S$ be a Seifert fibration as above. In case $M$ admits a spherical geometry, Kollár has constructed a uniruled real algebraic variety having a real component diffeomorphic to $M$. Therefore, we may assume that $M$ does

Both authors are members of the European Research Training Network RAAG (EC contract HPRN-CT-2001-00271). The first author is member of the European Research Training Network EAGER (EC contract HPRN-CT-2000-00099). The second author was supported in part by the CNRS, IM3, UMR 5030.

not admit a spherical geometry. Then, we show that there is a ramified Galois covering $p\colon \tilde{S} \to S$ of smooth surfaces such that the fiber product

$$\tilde{f}\colon \tilde{M} = M \times_S \tilde{S} \longrightarrow \tilde{S}$$

is a locally trivial circle bundle (Theorem 2.3). In particular, there is a finite group $G$ acting on the fiber bundle $\tilde{f}$ such that $\tilde{f}/G \cong f$. We show that there is a structure of a real algebraic surface on $\tilde{S}$, and that there is a real algebraic vector bundle $\tilde{L}$ of rank 2 on $\tilde{S}$ admitting

(1) a real algebraic action of $G$ on the total space of $\tilde{L}$, and
(2) a $G$-equivariant real algebraic metric $\lambda$ on $\tilde{L}$,

such that the unit circle bundle in $\tilde{L}$ is $G$-equivariantly diffeomorphic to $\tilde{M}$. The statement of Theorem 1.1 will then follow.

**Conventions.** A manifold is smooth, compact and connected, and without boundary, unless stated otherwise. A Riemann surface is compact and connected. A real algebraic variety is smooth projective and geometrically irreducible, unless stated otherwise.

**Acknowledgement.** We are grateful to János Kollár for pointing out a mistake in an earlier version of the paper. The second author wants to thank Jacques Lafontaine for helpful discussions.

## 2. Seifert fibrations

Let $S^1 \times D^2$ be the *solid torus* where $S^1$ is the unit circle $\{u \in \mathbb{C}, |u| = 1\}$ and $D^2$ is the closed unit disc $\{z \in \mathbb{C}, |z| \leq 1\}$. A *Seifert fibration* of the solid torus is a smooth map

$$f_{p,q} : S^1 \times D^2 \to D^2 \ , \ (u,z) \mapsto u^q z^p \ ,$$

where $p, q$ are relatively prime integers satisfying $0 \leq q < p$.

**Definition 2.1.** *Let $M$ be a 3-manifold. A* Seifert fibration *of $M$ is a smooth map $f : M \to S$ to a surface $S$ having the following property. Every $P \in S$ has a closed neighborhood $U$ such that the restriction of $f$ to $f^{-1}(U)$ is diffeomorphic to a Seifert fibration of the solid torus. More precisely, there are relatively prime integers $p$ and $q$, with $0 \leq q < p$, and there are diffeomorphisms $g\colon U \to D^2$ and $h\colon f^{-1}(U) \to S^1 \times D^2$ such that the diagram*

$$\begin{array}{ccc} f^{-1}(U) & \xrightarrow{h} & S^1 \times D^2 \\ \downarrow{\scriptstyle f_{|f^{-1}(U)}} & & \downarrow{\scriptstyle f_{p,q}} \\ U & \xrightarrow{g} & D^2 \end{array}$$

*commutes. We will say that $M$ is a* Seifert manifold *if $M$ admits a Seifert fibration.*

In the literature, e.g. [Sc83], nonorientable local models are also allowed. Following Kollár, we kept Seifert's original definition of a Seifert manifold.

Let us show that Seifert fibrations, as defined above, satisfy the property mentioned in the Introduction. The converse is easy to prove, and is left to the reader.

**Proposition 2.2.** *Let $f\colon M \to S$ be a Seifert fibration. Then, for every $P \in S$, there is a ramified covering $\tilde{U} \to U$ of an open neighborhood $U$ of $P$ such that the the fiber product $M \times_S \tilde{U} \to \tilde{U}$ is a locally trivial smooth circle bundle over $\tilde{U}$.*

Note that, the fiber product $M \times_S \tilde{U}$ is to be taken in the category of smooth manifolds. We will see that, in general, it does not coincide with the set-theoretic fiber product!

*Proof.* According to the definition of a Seifert fibration, it suffices to show that there is a ramified covering $g\colon D^2 \to D^2$ such that the fiber product

$$\tilde{f}\colon \tilde{F} = (S^1 \times D^2) \times_{f_{p,q},D^2,g} D^2 \longrightarrow D^2$$

is a trivial smooth circle bundle over $D^2$.

Let $g\colon D^2 \to D^2$ be the ramified covering $g(w) = w^p$. Then, the set theoretic fiber product of $S^1 \times D^2$ and $D^2$ over $D^2$ is the set

$$F = \{(u,z,w) \in S^1 \times D^2 \times D^2 | u^q z^p = w^p\}.$$

It is easy to see that $F$ is not necessarily smooth along the subset $S^1 \times \{0\} \times \{0\}$. Let $f\colon F \to D^2$ be the map defined by $f(u,z,w) = w$. Let $g'\colon F \to S^1 \times D^2$ be the map defined by $g(u,z,w) = (u,z)$, so that the diagram

$$\begin{array}{ccc} F & \xrightarrow{g'} & S^1 \times D^2 \\ \downarrow f & & \downarrow f_{p,q} \\ D^2 & \xrightarrow{g} & D^2 \end{array}$$

commutes.

The fiber product $\tilde{F}$ is the manifold defined by

$$\tilde{F} = \{(u,z,x) \in S^1 \times D^2 \times S^1 | u^q = x^p\}.$$

The desingularization map is the map $\delta\colon \tilde{F} \to F$ defined by $\delta(u,z,x) = (u,z,xz)$. Let $\tilde{f}\colon \tilde{F} \to D^2$ be the map defined by $\tilde{f}(u,z,x) = xz$. Then, the diagram

$$\begin{array}{ccc} \tilde{F} & \xrightarrow{\delta} & F \\ \downarrow \tilde{f} & & \downarrow f \\ D^2 & \xrightarrow{id} & D^2 \end{array}$$

commutes, i.e., the fibration $\tilde{f}$ is the desingularization of $f$.

Now, the map $\tilde{f}$ is a trivial fibration. Indeed, $\tilde{F}$ is a homogeneous space over $D^2$ under the action of $S^1$ defined by $v \cdot (u,z,x) = (v^p u, v^{-q} z, v^q x)$. Moreover, $\tilde{f}$ admits a smooth section $s$ defined by $s(z) = (1, z, 1)$. Therefore, $\tilde{f}$ is a trivial fibration of $\tilde{F}$ over $D^2$. $\square$

A point $P$ on a 2-dimensional orbifold $S$ is a *cone point* with *cone angle* $2\pi/p$ if a neighbourhood of $P$ is orbifold diffeomorphic to the orbifold quotient $\mathbb{C}//\mu_p$, where $\mu_p$ is the group of $p$-th roots of unity. A cone point is a *trivial* cone point, or a *smooth* cone point, if its cone angle is equal to $2\pi$.

Let $f\colon M \to S$ be a Seifert fibration. It follows from the local description of $f$ that the surface $S$ has a natural structure of an orbifold with only finitely many nontrivial cone points. Indeed, with the notation above, the fibration $f_{p,q}$ is the quotient of the trivial fibration $\tilde{f}$ by the action of $\mu_p$ defined on $\tilde{F}$ by $\xi \cdot (u,z,x) = (u,z,\xi x)$. The target of $f_{p,q}$ acquires the orbifold structure of $D^2//\mu_p$.

Recall that a manifold $M$ *admits a geometric structure* if $M$ admits a complete, locally homogeneous metric. In that case, the universal covering space $M'$ of $M$ admits a complete homogeneous metric. The manifold $M$ has then a geometric structure *modeled* on the $(Isom(M'), M')$-geometry.

More generally, a *geometry* is a pair $(I, V)$ where $V$ is a simply connected manifold and $I$ a real Lie group acting smoothly and transitively on $V$ with compact point stabilisers. We will only consider geometries $(I, V)$ that *admit a compact quotient*, i.e., there is a subgroup $H \subset I$ such that the projection $V \to V/H$ is a covering map onto a compact quotient. Two geometries $(I, V)$ and $(I', V')$ are equivalent if $I$ is isomorphic to $I'$ and there is a diffeomorphism $\varphi\colon V\to V'$ which transform the action $\rho\colon I\to\mathrm{Diff}(V)$ onto the action $\rho'\colon I'\to\mathrm{Diff}(V')$. We restrict ourselves to geometries where $I$ is maximal. Namely, if $I'$ is a strict subgroup of $I$ and $\rho'\colon I'\to\mathrm{Diff}(V)$ is the restriction of $\rho\colon I\to\mathrm{Diff}(V)$, we will only consider the $(I,V)$ geometry.

Thurston has classified the 3-dimensional geometries: there are eight of them. In general, a 3-manifold does not possess a geometric structure. However, it turns out that every Seifert manifold admits a geometric structure and that the geometry involved is unique [Sc83, Sec. 4]. Let $M$ be a Seifert manifold, the geometry of $M$ is modeled on one of the six following models (see [Sc83] for a detailed description of each geometry):

$$S^3, S^2\times\mathbb{R}, E^3, \mathrm{Nil}, H^2\times\mathbb{R}, \widetilde{\mathrm{SL}_2\,\mathbb{R}}$$

where $E^3$ is the 3-dimensional euclidean space and $H^2$ is the hyperbolic plane.

The appropiate geometry for a Seifert bundle can be determined from the two invariants $\chi$ and $e$, where $\chi$ is the Euler number of the base orbifold and $e$ is the Euler number of the Seifert bundle [Sc83, Table 4.1].

|         | $\chi > 0$       | $\chi = 0$ | $\chi < 0$                     |
|---------|------------------|------------|--------------------------------|
| $e = 0$ | $S^2\times\mathbb{R}$ | $E^3$      | $H^2\times\mathbb{R}$          |
| $e\neq 0$ | $S^3$           | Nil        | $\widetilde{\mathrm{SL}_2\,\mathbb{R}}$ |

**Table 1:** Geometries for Seifert manifolds.

**Theorem 2.3.** *Let $M$ be an orientable Seifert manifold that does not admit a spherical geometry, and let $f\colon M\to S$ be a Seifert fibration. Then there is an orientable surface $\tilde{S}$ and a finite ramified Galois covering $\tilde{S}\to S$ such that the fiber product*

$$\tilde{f}\colon \tilde{M} = M\times_S \tilde{S} \longrightarrow \tilde{S}$$

*is a locally trivial smooth circle bundle.*

A group of isometries of a Riemannian manifold $B$ is *discrete* if for any $x\in B$, the orbit of $x$ intersects a small neighborhood of $x$ only finitely many times. The quotient $B/\Gamma$ of $B$ by a discrete group $\Gamma$ of isometries, is a surface. The projection map $B\to B/\Gamma$ is a local homeomorphism except at points $x$ where the isotropy subgroup $\Gamma_x$ is nontrivial. In that case, $\Gamma_x$ is a cylic group $\mathbb{Z}/p\mathbb{Z}$ for some $p>1$, and the projection is similar to the projection of a meridian disk cutting across a singular fiber of a Seifert fibration.

For convenience, let us denote by $S^2(p,q)$, $1\leq q<p$, the orbifold whose underlying surface is $S^2$ with two cone points with angles $2\pi/q$ and $2\pi/p$. If $q=1$, $S^2(p,1)=S^2(p)$ is the teardrop orbifold. From [Sc83, Th. 2.3, Th. 2.4 and Th. 2.5], we can state the following:

**Theorem 2.4** (Scott). *Every closed $2$-dimensional orbifold $S$ with only cone points, and which is different from $S^2(p,q)$, $p\neq q$, is finitely covered by a smooth surface $\tilde{S}$.*

For the convenience of the reader we recall the main ideas of the proof.

*Proof.* Every 2-dimensional orbifold with only cone points, and which is different from $S^2(p,q)$, $p\neq q$, is isomorphic, as an orbifold, to the quotient of $S^2$, $E^2$ or $H^2$ by some discrete group of isometries $\Gamma$. We need to show that any finitely generated, discrete group $\Gamma$ of isometries of $S^2$, $E^2$ or $H^2$ with compact quotient space contains a torsion free subgroup of finite index. This is trivial for $S^2$ and easy for $E^2$. For $H^2$ it is a corollary of Selberg's Lemma below.

Let us denote by $\Gamma' \subset \Gamma$ a torsion-free normal subgroup of finite index, the orbifold quotient $\tilde{S} = S^2/\Gamma'$, $E^2/\Gamma'$ or $H^2/\Gamma'$, respectively, is then a smooth surface. □

**Selberg's Lemma.** [Ra94, Chap. 7] *Every finitely generated subgroup of $\mathrm{GL}_n(\mathbb{C})$ contains a torsion-free normal subgroup of finite index.* □

*Remark* 2.5. The fact that any finitely generated, discrete group of isometries of $H^2$ admits a torsion free subgroup of finite index was conjectured by Fenchel and the first proof was completed by Fox [Fo52]. From the modern point of view, however, this result is a corollary of the Selberg's Lemma

*Proof of Theorem 2.3.* Let $M$ be an orientable Seifert manifold that does not admit a spherical geometry. Let $M \to S$ be a Seifert fibration. Let us show that the base orbifold $S$ is not isomorphic to one of the orbifolds $S^2(p,q)$, with $1 \leq q < p$. Indeed, the Euler number $\chi$ of $S^2(p,q)$ is strictly positive, and the Euler number $e$ of any Seifert fibration over $S^2(p,q)$ is nonzero (see [Sc83], in particular, the discussion before Lemma 3.7). Therefore, by Table 1, if $S \cong S^2(p,q)$, the manifold $M$ would admit a spherical geometry. This shows that $S$ is not isomorphic to $S^2(p,q)$. It follows from Theorem 2.4 that there is a finite ramified Galois covering $\tilde{S} \to S$ of the orbifold $S$ by a smooth surface $\tilde{S}$. Moreover, we may assume $\tilde{S}$ to be orientable, taking the Galois closure of the orientation double covering if necessary. It is clear that the fiber product $\tilde{f}$ is locally trivial. □

## 3. Klein surfaces

In this section we recall the definition of a Klein surface, and we prove some statements that we need for the proof of Theorem 1.1. Classically, a Klein surface is defined as a topological surface endowed with an atlas whose transition functions are either holomorphic or antiholomorphic [AG71]. This seems to be a less suitable point of view for what we need since, with that definition, a Klein surface is not a locally ringed space. In particular, the definition of a line bundle over such a Klein surface is cumbersome, and, we would not have at our disposal a first cohomology group of the type $H^1(S, \mathcal{O}^\star)$ classifying all line bundles on a given Klein surface $S$. Therefore, we will use another definition of a Klein surface, giving rise to a category equivalent to the category of Klein surfaces of [AG71].

Let $\mathbb{D}$ be the double open half plane $\mathbb{C}\setminus\mathbb{R}$. On $\mathbb{D}$ one has the sheaf of holomorphic functions $\mathcal{H}$. The Galois group $\Sigma = \mathrm{Gal}(\mathbb{C}/\mathbb{R})$ acts naturally on $\mathbb{D}$. Let $\sigma$ denote complex conjugation in $\Sigma$. We consider the following *algebraic* action of $\Sigma$ on the sheaf $\mathcal{H}$ over the action of $\Sigma$ on $\mathbb{D}$. If $U \subseteq \mathbb{D}$ is open and $f$ is a section of $\mathcal{H}$ over $U$, then we define $\sigma \cdot f \in \mathcal{H}(\sigma \cdot U)$ by

$$(\sigma \cdot f)(z) = \overline{f(\overline{z})}$$

for all $z \in \sigma \cdot U$. Let $(\mathbb{H}, \mathcal{O})$ be the quotient of $(\mathbb{D}, \mathcal{H})$ by the action of $\Sigma$ in the category of locally ringed spaces. In particular, $\mathbb{H} = \mathbb{D}/\Sigma$ is homeomorphic to the open upper half plane—or lower half plane for that matter. Let $p\colon \mathbb{D} \to \mathbb{H}$ be the quotient map. Then $\mathcal{O}$ is the sheaf $(p_\star \mathcal{H})^\Sigma$ of $\Sigma$-invariant sections over $\mathbb{H}$. The sheaf $\mathcal{O}$ on $\mathbb{H}$ is a sheaf of local $\mathbb{R}$-algebras. Each stalk of $\mathcal{O}$ is noncanonically isomorphic to the $\mathbb{R}$-algebra $\mathbb{C}\{z\}$ of complex convergent power series in $z$.

A *Klein surface* is a locally ringed space $(S, \mathcal{O})$, where $\mathcal{O}$ is a sheaf of local $\mathbb{R}$-algebras, such that $S$ is compact connected and separated, and $(S, \mathcal{O})$ is locally isomorphic to $(\mathbb{H}, \mathcal{O})$. With the obvious definition of morphisms of Klein surfaces, we have the category of Klein surfaces. Note that, here, we have only defined the notion of a compact connected Klein surface without boundary. For a more general definition of Klein surfaces, the reader may refer to [Hu02].

Basic examples of Klein surfaces are the following. A Riemann surface is a Klein surface. More generally, let $S$ be a Riemann surface. A *Klein action* of a finite group $G$ on $S$ is an action of $G$ on $S$ as a Klein surface. Let be given a Klein action of $G$ on $S$. Suppose that $S$ contains only finitely many fixed points for the action of $G$. Then the quotient $S/G$ has a natural structure of a Klein surface.

The following statement is well known.

**Theorem 3.1.** *Let $S$ be a compact connected smooth surface. Then $S$ admits the structure of a Klein surface.*

*Proof.* If $S$ is orientable, then $S$ admits the structure of a Riemann surface. In particular, $S$ admits the structure of a Klein surface. Therefore, we may assume that $S$ is not orientable. This means that $S$ is a $(g+1)$-fold connected sum of real projective planes, for some natural integer $g$. Let $C$ be any real algebraic curve of genus $g$ without real points. Then, the set of complex points $C(\mathbb{C})$ of $C$ is a Riemann surface of genus $g$. Since $C$ is real, the Galois group $\Sigma = \mathrm{Gal}(\mathbb{C}/\mathbb{R})$ acts on $C(\mathbb{C})$ by holomorphic or antiholomorphic automorphisms. Since $C$ has no real points, the action of $\Sigma$ on $C(\mathbb{C})$ is fixed point-free, i.e., we are in the presence of a Klein action of $\Sigma$ on $S$. Therefore, the quotient $C(\mathbb{C})/\Sigma$ has a natural structure of a Klein surface. It is clear that $C(\mathbb{C})/\Sigma$ is diffeomorphic to $S$ as a smooth surface. Hence, $S$ admits the structure of a Klein surface. $\square$

There is also a Klein version of the Riemann Existence Theorem. It can either be proven as the Riemann Existence Theorem, or it can be proven using the Riemann Existence Theorem.

**Theorem 3.2.** *Let $S$ be a Klein surface and let $\tilde{S}$ be a compact connected smooth surface. If $f\colon \tilde{S} \to S$ is a ramified covering of smooth surfaces, then there is a unique structure of a Klein surface on $\tilde{S}$ such that $f$ is a morphism of Klein surfaces.* $\square$

Let $(S, \mathcal{O})$ be a Klein surface. A *line bundle* over $S$ is an invertible sheaf of $\mathcal{O}$-modules. The group of isomorphism classes of line bundles is isomorphic to the group $H^1(S, \mathcal{O}^\star)$.

We will also need the notion of a smooth line bundle over $S$. Let $\mathcal{C}^\infty$ be the sheaf of smooth complex valued functions on the open double half plane $\mathbb{D}$. The Galois group $\Sigma = \mathrm{Gal}(\mathbb{C}/\mathbb{R})$ acts on $\mathcal{C}^\infty$, in a similar way as its action on $\mathcal{H}$. This action extends the action of $\Sigma$ on $\mathcal{H}$. Denote by $\mathcal{C}$ the induced sheaf of $\Sigma$-invariant sections on $\mathbb{H}$. We call it the sheaf of *smooth functions* on $\mathbb{H}$. The sheaf $\mathcal{C}$ on $\mathbb{H}$ contains $\mathcal{O}$ as a subsheaf. It is now clear that a Klein surface $(S, \mathcal{O})$ carries an induced sheaf $\mathcal{C}$ of smooth functions which contains the sheaf $\mathcal{O}$ as a subsheaf.

A *smooth line bundle* on a Klein surface $(S, \mathcal{O})$ is an invertible sheaf of $\mathcal{C}$-modules. Again, the group of isomorphism classes of smooth line bundles on $S$ is isomorphic to $H^1(S, \mathcal{C}^\star)$. Of course, if $L$ is a line bundle on $S$, then $L \otimes_\mathcal{O} \mathcal{C}$ is a smooth line bundle on $S$. Let $L'$ be a smooth line bundle on $S$. We say that $L'$ admits the structure of a *Klein bundle* if there is a line bundle $L$ over $S$ such that $L \otimes_\mathcal{O} \mathcal{C} \cong L'$. The following statement shows that every smooth line bundle on a Klein surface does admit the structure of a Klein bundle.

**Theorem 3.3.** *Let $(S, \mathcal{O})$ be a Klein surface and let $\mathcal{C}$ be the induced sheaf of smooth functions on $S$. If $L'$ is a smooth line bundle on $S$ then there is a line bundle $L$ on $S$ such that*
$$L \otimes_\mathcal{O} \mathcal{C} \cong L'.$$

*Proof.* We show that the natural map
$$H^1(S, \mathcal{O}^\star) \longrightarrow H^1(S, \mathcal{C}^\star)$$

is surjective. As in [Hu02], we have exponential morphisms
$$\exp\colon \mathcal{O} \to \mathcal{O}^\star \quad \text{and} \quad \exp\colon \mathcal{C} \to \mathcal{C}^\star.$$
They are both surjective, and their kernels are isomorphic. Let $\mathcal{K}$ denote their kernels. Then we have a morphism of short exact sequences
$$\begin{array}{ccccccccc} 0 & \to & \mathcal{K} & \to & \mathcal{O} & \to & \mathcal{O}^\star & \to & 0 \\ & & \downarrow & & \downarrow & & \downarrow & & \\ 0 & \to & \mathcal{K} & \to & \mathcal{C} & \to & \mathcal{C}^\star & \to & 0 \end{array}$$
It induces the following commutative diagram with exact rows.
$$\begin{array}{ccccccc} H^1(S,\mathcal{O}) & \to & H^1(S,\mathcal{O}^\star) & \to & H^2(S,\mathcal{K}) & \to & H^2(S,\mathcal{O}) \\ \downarrow & & \downarrow & & \downarrow & & \downarrow \\ H^1(S,\mathcal{C}) & \to & H^1(S,\mathcal{C}^\star) & \to & H^2(S,\mathcal{K}) & \to & H^2(S,\mathcal{C}) \end{array}$$
Now, $\mathcal{C}$ is a fine sheaf. Hence, $H^1(S,\mathcal{C}) = 0$, and the map
$$H^1(S,\mathcal{C}^\star) \longrightarrow H^2(S,\mathcal{K})$$
is injective. Moreover, $H^2(S,\mathcal{O}) = 0$ [Hu02]. Hence, the map
$$H^1(S,\mathcal{O}^\star) \longrightarrow H^2(S,\mathcal{K})$$
is surjective. It follows that the natural map
$$H^1(S,\mathcal{O}^\star) \longrightarrow H^1(S,\mathcal{C}^\star)$$
is surjective. □

## 4. Equivariant line bundles on Riemann surfaces

Let $S$ be a Riemann surface and let $L$ be a smooth complex line bundle on $S$. Let be given a Klein action of a finite group $G$ on $S$. A *smooth Klein action* of $G$ on $L$ is an action of $G$ on $L$ over the action of $G$ on $S$ such that $g \in G$ acts antilinearly on $L$ if and only if $g$ acts antiholomorphically on $S$, for all $g \in G$. If, moreover, $L$ is a holomorphic line bundle on $S$ and $G$ acts by holomorphic or antiholomorphic automorphisms on the total space $L$, then the smooth Klein action is a *Klein action* of $G$ on $L$.

**Theorem 4.1.** *Let $S$ be a Riemann surface and let $L$ be a smooth complex line bundle over $S$. Let be given a faithful Klein action of a finite group $G$ on $S$ and a smooth Klein action of $G$ on $L$. Then, there is a structure of a holomorphic line bundle on $L$ such that the smooth Klein action of $G$ is a Klein action of $G$ on $L$.*

*Proof.* Since the action of $G$ on $S$ is a Klein action, $S$ contains finitely many fixed points $P_1, \ldots, P_n$. Let $G_1, \ldots, G_n$ be the isotropy groups of $P_1, \ldots, P_n$, respectively. Since the action of $G$ on $S$ is a Klein action, each isotropy group $G_i$ is a finite cyclic group of order $p_i$, acting holomorphically on $S$. Let $\rho_i$ be the induced 1-dimensional representation of $G_i$ on the complex tangent space $T_iS$ of $S$ at $P_i$. Since $\rho_i$ is a faithful representation and since the induced action of $G_i$ on the fiber $L_i$ over $P_i$ is complex linear, there is a unique integer $q_i \in \{0, \ldots, p_i - 1\}$ such that the 1-dimensional representation $L_i$ is isomorphic to $\rho_i^{q_i}$.

Let $K$ be the holomorphic line bundle $\mathcal{O}(\sum q_i P_i)$ on $S$. It is clear that $K$ comes along with a Klein action of $G$. Then, $K \otimes L$ is a smooth complex line bundle on $S$ with a smooth Klein action of $G$ such that the group $G_i$ acts trivially on the fiber $K_i \otimes L_i$ over $P_i$. Now, it suffices to show that $K \otimes L$ admits a structure of a holomorphic complex line bundle such that the smooth action of $G$ on $K \otimes L$ is a Klein action. Therefore, replacing $L$ by $K \otimes L$, we may assume that $q_i = 0$, for $i = 1, \ldots, n$. More precisely, we may assume that, for each $i$, the action of $G_i$ on the fiber $L_i$ of $L$ over $P_i$ is trivial.

Let $S'$ be the quotient Klein surface $S/G$, and let $p\colon S \to S'$ be the quotient map. Since the action of $G_i$ on $L_i$ is trivial, there is a smooth line bundle $L'$ on $S'$ such that $p^\star L'$ is $G$-equivariantly isomorphic to $L$. By Theorem 3.3, $L'$ has a structure of a Klein bundle over $S'$. It follows that $p^\star L'$ has the structure of a holomorphic line bundle such that the action of $G$ on $p^\star L'$ is a Klein action. $\square$

## 5. Uniruled algebraic models

*Proof of Theorem 1.1.* Let $M$ be an orientable Seifert manifold. We show that there is a uniruled real algebraic variety $X$ such that $M$ is diffeomorphic to a connected component of $X(\mathbb{R})$. That statement is known to be true if $M$ admits a spherical geometry [Ko99, Ex. 10.4]. Therefore, we may assume that $M$ does not admit a spherical geometry.

Choose a Seifert fibration $f\colon M \to S$ of $M$. By Theorem 2.3, there is a ramified Galois covering
$$p\colon \tilde{S} \longrightarrow S$$
such that the Seifert fibration
$$\tilde{f}\colon \tilde{M} \longrightarrow \tilde{S},$$
obtained from $f$ by base change, is a locally trivial circle fibration. Moreover, we may assume $\tilde{S}$ and $\tilde{M}$ to be oriented. Let $G$ be the Galois group of $\tilde{S}$ over $S$. Then $G$ acts naturally on $\tilde{M}$ and $\tilde{f}$ is $G$-equivariant. The quotient of $\tilde{f}\colon \tilde{M} \to \tilde{S}$ by $G$ is isomorphic to $f\colon M \to S$ as Seifert fibrations. This means that there is a quotient map
$$\tilde{p}\colon \tilde{M} \longrightarrow M$$
for the action of $G$ on $\tilde{M}$ such that the diagram
$$\begin{array}{ccc} \tilde{M} & \xrightarrow{\tilde{p}} & M \\ \downarrow{\tilde{f}} & & \downarrow{f} \\ \tilde{S} & \xrightarrow{p} & S \end{array}$$
commutes.

Since $S$ is a compact connected surface without boundary, $S$ admits a structure of a Klein surface by Theorem 3.1. By the Riemann Existence Theorem for Klein surfaces, there is a unique structure of a Klein surface on $\tilde{S}$ such that the map $p\colon \tilde{S} \to S$ is a morphism of Klein surfaces. In particular, the group $G$ acts on $\tilde{S}$ by automorphisms of $\tilde{S}$. In fact, since $\tilde{S}$ is an oriented Klein surface without boundary, $\tilde{S}$ is a Riemann surface. The action of $G$ on the Riemann surface $\tilde{S}$ is by holomorphic or antiholomorphic automorphisms, i.e., we are in the presence of a so-called Klein action of $G$ on the Riemann surface $\tilde{S}$.

Choose a smooth relative Riemannian metric $\mu$ on $\tilde{M}/\tilde{S}$. Since $G$ is finite, one may assume that $\mu$ is $G$-equivariant. Since $\tilde{f}$ is locally trivial, and since $\tilde{M}$ and $\tilde{S}$ are oriented, there is a relative orientation of $\tilde{M}/\tilde{S}$. Hence, the structure group of the locally trivial circle bundle $\tilde{M}/\tilde{S}$ is $\mathrm{SO}(2)$. Since $\mathrm{SO}(2) = \mathrm{SU}(1)$, there is a smooth complex line bundle $\tilde{L}$ on $\tilde{S}$, that comes along with a hermitian metric, whose unit circle bundle is $\tilde{M}$. We also have an action of $G$ on $\tilde{L}$ over the action of $G$ on $\tilde{S}$ that extends the action of $G$ on $\tilde{M}$. The action of $G$ on $\tilde{L}$ is a smooth Klein action since $G$ acts by orientation preserving automorphisms on $\tilde{M}$. By Theorem 4.1, there is a structure of a homolomorphic line bundle on $\tilde{L}$ such that the action of $G$ is a Klein action. By the GAGA-principle, $\tilde{S}$ is a complex algebraic curve and $\tilde{L}$ is a complex algebraic line bundle on $\tilde{S}$. Moreover, the action of $G$ on $\tilde{L}$ is by algebraic or antialgebraic automorphisms. The restriction of scalars $R(\tilde{S})$ is a real algebraic surface whose set of real points is diffeomorphic to $\tilde{S}$ [Hu92, Hu00]. The restriction of scalars $R(\tilde{L})$ is a real algebraic vector bundle over $R(\tilde{S})$ of rank 2, whose set of real points is diffeomorphic to $\tilde{L}$. The action of $G$ on $\tilde{L}$ induces an algebraic action of $G$ on $R(\tilde{L})$.

Let $U$ be an affine Zariski open subset of $R(\tilde{S})$ containing the real points of $R(\tilde{S})$ and which is $G$-equivariant. Since $R(\tilde{L})$ is a real algebraic vector bundle over an affine real algebraic variety, there is a vector bundle $V$ over $U$ such that the direct sum
$$V \oplus (R(\tilde{L})_{|U})$$
is trivial. Since $V \oplus (R(\tilde{L})_{|U})$ is trivial, there is a real algebraic metric $\lambda$ on the restriction of $R(\tilde{L})$ to $U$. Since $G$ is finite, we may assume that $\lambda$ is $G$-equivariant. Let $\tilde{T}$ be the unit circle bundle $\lambda = 1$ in $R(\tilde{L})_{|U}$. Then $\tilde{T}(\mathbb{R})$ is $G$-equivariantly diffeomorphic to $\tilde{M}$. In particular, the quotient $\tilde{T}(\mathbb{R})/G$ is diffeomorphic to $M$. Let $X$ be a smooth projective model of $\tilde{T}/G$. Since $G$ acts fixed point-freely on $\tilde{T}(\mathbb{R})$, the quotient $\tilde{T}(\mathbb{R})/G$ is smooth and is, therefore, a connected component of $X(\mathbb{R})$. The real algebraic variety $X$ is uniruled since a smooth projective model of $\tilde{T}$ is ruled. □

Johannes Huisman, Institut Mathématique de Rennes, Université de Rennes 1 Campus de Beaulieu, 35042 Rennes Cedex, France, Tel: +33 2 23 23 58 87, Fax: +33 2 23 23 67 90
  *E-mail address*: huisman@univ-rennes1.fr

Frédéric Mangolte, Laboratoire de Mathématiques, Université de Savoie, 73376 Le Bourget du Lac Cedex, France, Phone: +33 (0)4 79 75 86 60, Fax: +33 (0)4 79 75 81 42
  *E-mail address*: mangolte@univ-savoie.fr